\documentclass{amsart}
\usepackage[all]{xy}
\usepackage{a4wide}
\newtheorem{theorem}{Theorem}[section]
\newtheorem{corollary}[theorem]{Corollary}
\newtheorem{lemma}[theorem]{Lemma}
\newtheorem{example}[theorem]{Examples}
\newtheorem{definition}[theorem]{Definition}
\newtheorem{proposition}[theorem]{Proposition}

\newcommand{\NN}{\mathbb{N}}

\newcommand{\Tr}[2]{\mathrm{Tr}{\left( {#1}, {#2} \right)}}

\newcommand{\End}[2]{\mathrm{End}_{#1}(#2)}
\newcommand{\Hom}[3]{\mathrm{Hom}_{#1}(#2,#3)}
\newcommand{\Ker}[1]{\mathrm{Ker}(#1)}

\newcommand{\udim}[1]{\mathrm{udim}(#1)}

\newcommand{\AH}{{A \# H}}

\renewcommand{\Im}{\mathrm{Im}}

\begin{document}

\title{Irreducible actions and compressible modules}

\author{In\^es Borges}
\address{Instituto Superior De Contabilidade e Administra\c{c}\~ao de Coimbra, Quinta Agricola - Bencanta, 3040-316 Coimbra, Portugal}

\author{Christian Lomp}
\dedicatory{dedicated to the memory of John Dauns}
\address{Centro de Matemática de Universidade do Porto, Rua Campo Alegre 687, 4169-007 Porto, Portugal}

\thanks{The first author was supported by grant SFRH/PROTEC/49857/2009. The second author was partially supported by Centro de Matemática da Universidade do Porto (CMUP), financed by FCT (Portugal) through the programs POCTI (Programa Operacional Ciência, Tecnologia, Inovação) and POSI (Programa Operacional Sociedade da Informação), with national and European community structural funds.}



\begin{abstract}
Any finite set of linear operators on an algebra $A$ yields an operator algebra $B$ and a module structure on $A$, whose endomorphism ring is isomorphic to a subring $A^B$ of certain invariant elements of $A$.
We show that if $A$ is a critically compressible left $B$-module, then the dimension of its
self-injective hull $\widehat{A}$ over the ring of fractions of $A^B$ is bounded by the uniform dimension of $A$ and the number of linear operators
generating $B$. This extends a known result on irreducible Hopf actions and applies in particular to weak Hopf action.
Furthermore we prove necessary and sufficient conditions for an algebra $A$ to be critically compressible in the case of group actions, group gradings and
Lie actions.
 \end{abstract}

\maketitle


%
%


\renewcommand{\Im}{\mathrm{Im}}

\section{Introduction}
The starting point of this note is a theorem by Bergen et al. that says that for a Hopf algebra $H$ acting on an algebra $A$ with subring of invariants $A^H$, such that $A$ is an irreducible $\AH$-module, the dimension $dim(A_{A^H})$ is bound by the uniform dimension of $_AA$ and the number of generators of $\AH/\mathrm{Ann}_{\AH}(A)$ as $A$-module. We will generalize their result in two directions: replacing the Hopf action by a subalgebra of linear operators of $A$ which contains the regular left action of $A$ on itself and by weakening the irreducible assumption of $A$ to $A$ being a monoform respectively critically compressible module.

Actions of Hopf algebras include group actions, Lie actions and finite group gradings. Several generalizations of Hopf algebras
have emerged in recent years, like weak Hopf algebras (or quantum groupoids) introduced by B\"ohm et al. \cite{BohmNilSzlachanyi}. 
Examples of actions of weak Hopf algebras are given by groupoid actions on $C^*$-algebras (see \cite{Vallin}) and by quantum groupoids  arising from
Jones towers (see \cite{KadisonNikshych}). 
Our point of view is that such action on an algebra $A$ should be studied by looking at the module structure on $A$ given by certain algebras of linear operators, i.e. subalgebras of $\End{k}{A}$. Let $k$ be a commutative ring and let $A$ be an associative unital $k$-algebra. For any $a\in A$ define two linear operators $L_a$ and $R_a$ in $\End{k}{A}$ given by $L_a(x)= ax$ and $R_a(x) = xa$ for all $x\in A$. We identify $A$ with  the subalgebra $L(A)$ of $\End{k}{A}$ generated by all left multiplications $L_a$ and denote the subalgebra generated by all operators $L_a$ and $R_a$ by $M(A)$, which is also sometimes referred to as the {\it multiplication algebra} of $A$. As a left $L(A)$-module, $A$ is isomorphic to $L(A)$ since we assume $A$ to be unital. 

We will be interested in certain actions on an algebra $A$ that may stem from a bialgebra or more generally a bialgebroid. 
Usually we are interested in extensions $A\subseteq B$ where $B$ act on $A$ through a ring homomorphism $\phi:B\rightarrow
End_k (A)$ such that $\phi(a)=L_a$ for all $a \in A$ (like for instance if $B=\AH$ is a smash product). To study the intrinsic properties of $A$  under
this action it is enough to look at the subalgebra $\phi(B)$ in $\End{k}{A}$ generated by this action. Therefore we will consider subalgebras $B$ of $\End{k}{A}$ that contain $L(A)$ and act on $A$ by evaluation. Thus $A$ becomes a cyclic faithful left $B$-module by $\varphi\cdot a := \varphi(a)$, $\forall \varphi\in B,a\in A$. Note that $L_{a'}\cdot a = a'a$ for any $a,a'\in A$. 

Since we assume $A$ to be unital, the map $\Psi:\End{B}{A} \rightarrow A$ with $\Psi(f)=(1)f$, for all $f\in \End{B}{A}$, 
is an injective ring homomorphism, since $\forall f,g \in \End{B}{A}:$ 
$$ \Psi(f\circ g) = ((1)f)g = ((1)f\cdot 1)g = (1)f \cdot (1)g = \Psi(f)\Psi(g).$$
Note that we will write homomorphisms opposite to scalars.  
Moreover if $\Psi(f)=(1)f=0$, then $(a)f=a(1)f=0$ and $f=0$. The subalgebra $\Psi(\End{B}{A})$, which we will denote by $A^B$, can be described as the
set of elements $a\in A$ such that for any $b\in B: b\cdot a = (b\cdot 1)a$. 


A subset $I$ of $A$ is called $B$-stable if $B\cdot I \subseteq I$. The $B$-stable left ideals are precisely the (left) $B$-submodules
of $A$. In particular,  by restricting $\Psi$, we have $\Hom{B}{A}{I}\simeq I \cap A^B$, for any $B$-stable left ideal of $A$.

\begin{example}
Let $H$ be a weak Hopf algebra over $k$ (see Definition \ref{weakHopfDef}) and let $A$ be a left $H$-module algebra. Denote the
action of an element $h\in H$ on $A$ by $\lambda_h\in \End{k}{A}$ and define $B$ as the subalgebra of $\End{k}{A}$ generated by $L(A)$ and all operators $\lambda_h$.
The left $B$-submodules of $A$ are precisely the $H$-stable left ideals of $A$. Moreover $B$ is a factor of the smash product $A\# H$. A particular
case of this are  group actions given by a group $G$ 
and a group homomorphism $\eta: G \rightarrow Aut_k(A)$. Set $B=\langle L(A)\cup \eta(G) \rangle\subseteq \End{k}{A}$. The $B$-submodules of $A$ are precisely the $G$-stable left ideals of $A$ and $\Hom{B}{A}{I} \simeq I\cap A^G$, where
$A^G = \{ a\in A\mid  \forall g\in G: \eta(g)(a) = a \}$ is the fix ring of $A$. $B$ is a quotient of the skew group ring $A*G$.
Another example is given by Lie algebras acting as derivations on $A$, e.g. if $\delta \in Der_k(A)$ is a $k$-linear derivation of $A$ then consider $B=\langle L(A)\cup \{ \delta\}\rangle \subseteq \End{k}{A}$. The left $B$-submodules of $A$ are the left ideals $I$ that satisfy $\delta(I)\subseteq I$. The operator algebra $B$ is a factor of the  ring of differential operator $A[z;\delta]$, which as a left $A$-module is equal to $A[z]$ and its multiplication is given by $za=az+\delta(a)$.
The map $A[z;\delta] \rightarrow B$ with  $\sum_{i=0}^n a_iz^i \mapsto \sum_{i=0}^n L_{a_i}\circ \delta^i \in B$ is a surjective $k$-algebra
homomorphism and for any left $A[z;\delta]$-module $M$ we have $\Hom{A[z;\delta]}{A}{M}=\Hom{B}{A}{M}=M^\delta=\{m\in M\mid zm=0\}$. In particular
$\End{A[z;\delta]}{A}\simeq A^\delta = \Ker{\delta}$.
\end{example}

\begin{example} For $B=M(A)$, the $B$-submodules of $A$ are the two-sided ideals  of $A$, and  $\Hom{B}{A}{I}\simeq Z(A)\cap I$, where $Z(A)$ denotes
the centre of $A$. The algebra $M(A)$ is a quotient of the enveloping algebra $A^e=A\otimes A^{op}$ of $A$. 

Let $*$ be an involution of $A$. Set $B=\langle L(A)\cup\{*\}\rangle$. The left $B$-submodules of $A$ are the twosided $*$-ideals and $A^B=Z(A,*)$
is the subring of central symmetric elements of $A$. Note that $B$ can be seen as the factor ring of the skew-group ring $A^e \# G$ where $G=\{id,
\overline{*}\}$ is the cyclic group of order two and $\overline{*}\in Aut(A^e)$ is given by $ (a\otimes b)^{\overline{*}} := b^* \otimes a^*.$
\end{example}

\section{Irreducible actions of linear Operators.}
The left uniform dimension of $A$ is denoted by $\udim{A}$ and is the supremum of the cardinalities of the index sets of direct sums of left ideals
contained in $A$. A result by Bergen et al. \cite[Theorem 2.2]{BergenCohenFishman} says that if a Hopf algebra $H$ acts finitely on a
module algebra $A$ with finite uniform dimension, such that $A$ is a simple $\AH$-module, then $A$ has finite dimension over $A^H$.

The argument of \cite[Theorem 2.2]{BergenCohenFishman} uses the Jacobson Density Theorem which had been generalized by
Zelmanowitz in \cite{Zelmanowitz_GeneralDensity, Zelmanowitz}. The hypotheses of Zelmanowitz' density theorems are weaker than assuming the
existence of a  faithful simple module: a non-zero left $R$-module $M$ is called \textbf{compressible} if it can be embedded in each of its non-zero
submodules and it is called  \textbf{critically compressible} if it is compressible and cannot be embedded in any of its proper factor modules.
A left $R$-module $M$ is called {\bf monoform} if any non-zero partial endomorphism $f:N\rightarrow M$  from a submodule $N$ to $M$ is injective.
It is easy to see that the critically compressible modules are precisely those that are compressible and monoform.  Zelmanowitz proved the following two weak Density Theorem for rings which have a faithful critically compressible module respectively faithful monoform module.

\begin{theorem}[Zelmanowitz, {\cite[Theorem 2.2]{Zelmanowitz_WeaklyPrimitive}},{\cite[Proposition 2.1]{Zelmanowitz}}]\label{m-dense}\label{prop4}
Let $M$ be a faithful left $R$-module with self-injective hull $\widehat{M}$ and $\Delta=\mathrm{End}_R(\widehat{M})$.
\begin{enumerate}
  \item If $_RM$ is critically compressible, then for any elements $v_1,v_2,\ldots,v_k \in \widehat{M}$ that are lineary independent over $\Delta$ there exists $0 \neq f \in \Delta$ such that for any elements $u_1,u_2,\ldots,u_k \in M$ there exists $r \in R$ with $rv_i=(u_i)f$, for each $i=1,\ldots k$.
  \item If $_RM$ is monoform left $R$-module, then for any $v_1,\ldots,v_k \in \widehat{M}$, linearly independent over $\Delta$ and $u_1,\ldots,u_k \in \widehat{M}$ with
$u_1\neq 0$,  there exist $r , s \in R$ with $r u_j= s v_j$  for $j=1,\ldots, k$ and $r u_1 \neq 0$.
\end{enumerate}
A ring $R$ with a faithful critically compressible left $R$-module is called left {\it weakly primitive}.
\end{theorem}

The self-injective hull of a left $R$-module $M$ is the submodule $\widehat{M}=\Tr{M}{E(M)}$, where $E(M)$ is the injective
hull of $M$ and  $$\Tr{X}{Y}= \sum \{\Im{f} \mid f \in \Hom{R}{X}{Y}\}$$ is the trace of $X$ in $Y$. Hence any element of $\widehat{M}$ can be written as a finite sum $\sum_{i=1}^n (m_i)f_i$ of images of homomorphisms $f_i:M\rightarrow E(M)$.  If $M$ is a simple module, then $\widehat{M}=M$.

Using the self-injective hull of $_BA$, we are able of extending \cite[Theorem 2.2]{BergenCohenFishman}.

\begin{theorem}\label{Proposition1}\label{prop5}\label{FiniteDimensional}
Suppose $L(A)\subseteq B \subseteq End_{k}(A)$ such that $_A{B}$ is generated by $n$ elements. If one of the following conditions hold
\begin{enumerate}
 \item $A$ is a domain and $A$ is a monoform left $B$-module {\bf or }
 \item $A$ is a critically compressible left $B$-module {\bf or }
 \item $A$ is a simple left $B$-module
\end{enumerate}
 then $dim(\hat{A}_\Delta) \leq n \cdot udim(A)$, where $\hat{A}$ is the self-injective hull of $_BA$ and $\Delta=End_{B}(\hat{A})$. 
Note that in case (3) $\hat{A}=A$ and in case (2)  $\Delta$  is isomorphic to the division ring of fractions of $A^B$. 
\end{theorem}

\begin{proof}
We want to show that the dimension of $\widehat{A}_\Delta$ is bounded by $n \cdot udim(A)$. 
Since $\widehat{A}$ is generated as $\Delta$-vector space by the elements of $A$ 
(any element is a linear combination $\sum_{i=1}^n (a_i)f_i$ with $f_i\in \Delta$ and $a_i\in A$) 
it is enough to consider linearly independent elements in $A$. 
Let $v_1,\ldots,v_k \in A$ be linearly independent over $\Delta$. 
Under the condition (1) or (2), we will show that  for each $1\leq i\leq k$ there exist $s_i \in B$ and $a_i\in A$ such that 
$$s_i\cdot v_j = \delta_{ij} x_i,$$ where $x_i$ is a right non-zero divisor of $A$. Once the existence of such elements is guaranteed, we 
can proceed as follows: since $B$ is $n-$generated over $A$, there exist an epimorphism of left $A$-modules $\varphi:A^{n}\rightarrow B$. Hence for each $1\leq i \leq k$, there exist an element $t_i \in A^{n}$ with  $(t_i)\varphi = s_i$. We will show that $\sum_{i=1}^k At_i$ is direct. 
Suppose that there are elements $a_1,\ldots, a_k \in A$ such that $\sum_{i=1}^{k} a_i t_i =0$, then for all $1\leq j \leq k$: 
$$0=\left(\sum_{i=1}^{n} a_i t_i\right)\varphi \cdot v_j = \sum_{i=1}^{n} a_i \left(t_i\right)\varphi \cdot v_j =\sum_{i=1}^{n} a_i \left(s_i \cdot
v_j\right) =\sum_{i=1}^{n} a_i \delta_{ij} x_i = a_jx_j.$$
And since $x_j$ is a right non-zero divisor, $a_j=0$. Thus $\bigoplus_{i=1}^{n} A t_i$ is a direct sum in $A^n$ and the result follows since 
$$dim(\hat{A}_\Delta) \leq udim(A^n) = n \cdot udim(A).$$

To guarantee the existence of such elements $s_i$ and $x_i$ we use Zelmanowitz' weak density theorems:
Given $v_1,\ldots,v_k \in A$ be linearly independent over $\Delta$, define $u_{ij}= \delta_{ij}1_A$, for all $1\leq i,j\leq k$.

(1) Suppose that $A$ is a domain and $_BA$ is monoform. Since $_BA$ is faithful, Theorem \ref{m-dense} shows that for each $1\leq i \leq k$, there exist $s_i,r_i \in B$ with $$s_i\cdot v_j = r_i\cdot  u_{ij} = \delta_{ij} (r_i\cdot 1_A),$$ for all $1\leq j \leq k$, and $x_i=r_i\cdot 1_A \neq 0$. Since $A$ is a domain, $x_i$ is a right non-zero divisor.

(2) Suppose that $_BA$ is critically compressible, then Theorem \ref{prop4} shows that for each $1\leq i \leq k$, there exist 
$f_i\in \Delta$ and $s_i\in B$ such that $$ s_i\cdot v_j = (u_{ij})f_i = \delta_{ij} (1)f_i, \mbox{ for all } 1\leq j \leq k.$$ Since any $f\in
\Delta$ is injective, $a(1)f=(a)f\neq 0$, for any non-zero $a\in A$. Hence $x_i=(1)f_i$ is a right non-zero divisor in $A$.

By  \cite[11.5(2)]{Zelmanowitz_WeaklyPrimitive}, $\Delta$ is isomorphic to $\mathrm{Frac}(A^B)$.

(3) Follows from $(2)$.

\end{proof}

Theorem \ref{Proposition1} makes no statement if the uniform dimension of $A$ is infinite. However a domain with finite uniform
dimension has uniform dimension $1$ and is called left Ore domain (see \cite{GoodearlWarfield}). If $R$ is a
left Ore domain and $f:I\rightarrow R$ is a partial non-trivial left $R$-linear endomorphism, then there exists $x\in I$ such that $(x)f \neq 0$. For
any $y=rx\in \Ker{f}\cap Rx$ we have $r(x)f = (y)f=0$. Since $R$ is a domain and $(x)f\neq 0$, $y=0$ and $\Ker{f}\cap Rx=0$. Since $_RR$ is uniform
$\Ker{f}=0$, i.e. $f$ is injective. Hence any left Ore domain $R$ is itself a monoform left $R$-module. On the other hand if $R$ is a monform left
$R$-module, then it is also a  uniform module by \cite[11.3]{Wisbauer}. Moreover for any $a\in R$, the endomorphism $R_a:[x\mapsto xa]$ is injective. Hence $ba=0$
implies $b=0$ or $a=0$, i.e. $R$ is a domain. This shows the following Lemma:

\begin{lemma}\label{lemma1}
Let $A$ be $k$-algebra. Then $A$ is a left Ore domain if and only if $A$ is a monoform left $A$-module. In this case $A$ is a monoform left
$B$-module, for all $L(A)\subseteq B\subseteq \End{k}{A}$.
\end{lemma}

The last statement is clear, because if $f:I\rightarrow A$ is a $B$-linear partial endomorphism, then it is also $A$-linear and hence injective if
$A$ was monoform as left $A$-module. 

\begin{corollary}
 Let $A$ be a left Ore domain and $L(A)\subseteq B \subseteq \End{k}{A}$ with $_AB$ being generated by $n$ elements, then
$dim(\hat{A}_\Delta) \leq n$, where $\hat{A}$ is the self-injective hull of $A$ and $\Delta=End_{B}(\hat{A})$.
\end{corollary}

Before we apply this result to quantum groupoid actions, we first generalize another result from the theory of Hopf actions to a more general setting.

Let $A\subseteq B$ be an extension of $k$-algebras such that there exists a ring homomorphism $\Psi: B\rightarrow \End{k}{A}$ with $\Psi(a)=L_a$ for all $a\in A$. If $B'=\mathrm{Im}(\Psi)$, then $\End{B}{A}\simeq \Delta:=A^{B'}$. The extension $A\subseteq \AH$ with $H$ being a (weak) Hopf algebra acting on $A$ is an example of such an extension.
\begin{proposition}\label{Galois} Let $A\subseteq B$ be as above and assume that $B$ is a finitely generated left and right $A$-module. Let $\Delta=\End{B}{A}$. If $A$ has finite left Goldie dimension and $A$ is a simple left $B$-module, then the following statements are equivalent:
\begin{enumerate}
\item[(a)] $B$ is  a simple ring;
\item[(b)] $B\simeq \End{}{A_\Delta}$ and $A_\Delta$ is finitely generated projective;
\item[(c)] $A$ is a faithful left $B$-module and $\Hom{B}{A}{B} \neq 0$;
\end{enumerate}
If $B_A$ is free of rank $n$, then $(a-c)$ are also equivalent to:
\begin{enumerate}
\item[(d)]$dim(A_\Delta)=n$ and $\Hom{B}{A}{B} \neq 0$;
\end{enumerate}
\end{proposition}

\begin{proof} By Theorem \ref{FiniteDimensional}, $A$ has finite dimension over $\Delta$ and since $B$ is finitely generated as a left and right $A$-module, $B$ is an Artinian ring.

$(a)\Rightarrow (b)$ Suppose $B$ is simple, then it is  semisimple Artinian and the epimorphism $b\mapsto b\cdot 1_A$ splits as left $B$-module by some $B$-linear map $\phi:A\rightarrow B$, i.e. $_BA$ is projective. Since $B$ is simple and the trace ideal of $A$ in $B$ is non-zero, $\Tr{A}{B}=B$, i.e. $A$ is a generator in $B$-Mod. By a standard module theoretic argument $B\simeq \End{}{A_{\Delta}}$.

$(b)\Rightarrow (c)$ is clear since $A$ is a generator. 

$(c)\Rightarrow (a)$ $\Hom{B}{A}{B} \neq 0$ implies that $A$ is isomorphic to a minimal left ideal of $B$ (since $A$ is a simple $B$-module). 
As $_BA$ is faithful, $B$ is a left primitive ring having a faithful minimal left ideal, i.e. 
$B\simeq M_m(\Delta)$ with $m=dim(A_\Delta)$. 

If $B_A$ is free of rank $n$, then $n\cdot dim(A_\Delta) = dim(B_{\Delta})$ holds.

$(c)\Rightarrow (d)$ As seen in the last step $(c)$ implies $B\simeq M_m(\Delta)$ with $m=dim(A_\Delta)$. On the other hand 
$n\cdot m = dim(B_\Delta)$ holds, yields $n=m$. 

$(d)\Rightarrow (b)$ if $m=dim(A_\Delta)=n$, then $dim(B_\Delta)=m^2=dim(\End{}{A_\Delta})$. Thus $B\simeq \End{}{A_\Delta}$.
\end{proof}

We intend to apply Theorem \ref{FiniteDimensional} to weak Hopf algebra actions and recall therefore its definition  from \cite{BohmNilSzlachanyi}.

\begin{definition}\label{weakHopfDef}
An associative $k$-algebra $H$ with multiplication $m$ and unit $1$ which is also a coassociative coalgebra with comultiplication $\Delta$ and counit
$\epsilon$ is called a {\it quantum groupoid} or a {\it weak Hopf algebra} if it satisfies the following properties:
\begin{enumerate}
\item the comultiplication is multiplicative, i.e. for all $g,h \in H$: $\Delta(gh)=\Delta(g)\Delta(h)$.
\item the unit and counit satisfy $\epsilon(fgh)=\epsilon(fg_1)\epsilon(g_2h) = \epsilon(fg_2)\epsilon(g_1h)$ and 
\begin{equation*}(\Delta \otimes id)\Delta(1) = (\Delta(1)\otimes 1)(1\otimes \Delta(1)) = (1\otimes \Delta(1))(\Delta(1)\otimes 1)\end{equation*}
\item there exists a linear map $S:A\rightarrow A$, called {\it antipode}, such that for all $h\in H$: $$S(h)=S(h_1)h_2S(h_3),$$ 
$$h_1S(h_2)=(\epsilon \otimes id)(\Delta(1)(h\otimes 1)) =: \epsilon_t(h),$$
$$S(h_1)h_2=(id \otimes \epsilon)((1\otimes h)\Delta(1)) =: \epsilon_s(h).$$
\end{enumerate}
Note that we will use Sweedler's notation for the comultiplication with suppressed summation symbol.
\end{definition}

The image of $\epsilon_t$ and $\epsilon_s$ are subalgebras $H_t$ and $H_s$ of $H$ which are separable over $k$ (\cite[2.3.4]{NikshychVainerman})§. Those subalgebras are also characterized by $ H_t = \{h\in H : \Delta(h)= 1_1h\otimes 1_2\}$ respectively $ H_s
= \{h\in H : \Delta(h)= 1_1\otimes 1_2h\}$. A left $H$-module algebra $A$ over a quantum groupoid $H$ is an associative unital algebra $A$ such that
$A$ is a left $H$-module and for all $a,b\in A, h\in H$:
\begin{equation}\label{intertwining}h\cdot(ab)=(h_1\cdot a)(h_2\cdot b) \: \mbox{ and } \: h\cdot 1_A = \epsilon_t(h)\cdot 1_A\end{equation}
Let $A$ be a left $H$-module algebra over a quantum groupoid $H$ and let $\lambda$ be the ring homomorphism from $H$ to $\End{k}{A}$ that defines the
left module structure on $A$, i.e. $\lambda_h(a) := h\cdot a$ for all $h\in H, a\in A$. Property (\ref{intertwining}) of the definition
above can be interpreted as an intertwining relation $ \lambda_h\circ L_a = L_{h_1\cdot a}\circ \lambda_{h_2}$ of left multiplications $L_a$ and left
$H$-actions $\lambda_h$. Recall that the smash product $A\# H$ of a left $H$-module algebra $A$ and a quantum groupoid $H$ is defined on the tensor
product $A\otimes_{H_t} H$ where $A$ is considered  a right $H_t$-module by $a\cdot z = a(z\cdot 1_A)$ for $a\in A, z\in H_t$. Consider the algebra extension $A\subseteq \AH=:B$ and the subalgebra of $\End{k}{A}$
generated by the $H$-action and $L(A)$, we can apply Theorem \ref{FiniteDimensional} and Proposition \ref{Galois} to obtain an quantum groupoid analog of
results by Bergen et al. \cite[Theorem 2.2]{BergenCohenFishman} and Cohen et al. \cite[Theorem 3.3]{CohenFishmanMontgomery}. Recall, from \cite{CaenepeelDeGroot} that a $H$-comodule algebra $A$ is an $H$-Galois extension of $A^{coH}$ if the canonical map $A\otimes_{A^{coH}} A \rightarrow A\otimes H$ is an isomorphism. If $H$ is finite dimensional and $A$ a left $H$-module algebra such that $A^H$ is divison ring, then by \cite[Proposition 2.3]{CaenepeelDeGroot} $A$ is a projective generator in $\AH$-Mod if and only if $A^H\subseteq A$ is $H^*$-Galois.

\begin{corollary} Let $A$ be a left $H$-module algebra over a finite dimensional quantum groupoid $H$. 
\begin{enumerate}
\item If $A$ is a domain and monoform as left $\AH$-module or if $A$ is critically compressible as left $\AH$-module, then  $dim(A_{A^H})\leq
\mathrm{dim}(H)\cdot \udim{A}$.
\item  If $A$ has finite left uniform dimension and is
a simple left $A\# H$-module, then the following statements are equivalent:
\begin{enumerate}
\item[(a)] $\AH$ is simple;
\item[(b)] $A^H\subseteq A$ is an $H^*$-Galois extension;
\item[(c)] $A$ is a faithful left $\AH$-module;
\end{enumerate}
If $A\#H_A$ is free of rank $dim(H)$, then $(a-c)$ are also equivalent to:
\begin{enumerate}
\item[(d)] $dim(A_{A^H}) = dim(H)$.
\end{enumerate}
\end{enumerate}
\end{corollary}

The associativity of $A$ is not needed to prove Theorem \ref{FiniteDimensional}. Non-associative examples of module algebras are given by module
algebras over quasi-Hopf algebras. Let $H$ be a quasi-Hopf algebra, that is $H$ is an associative algebra which is a not necessarily coassociative
coalgebra satisfying some compatibility conditions (see \cite{Drinfeld}). $H$ acts on an algebra $A$ if $A$ is a unital algebra in the category of left
$H$-modules (see \cite{Bulacu}). In particular its multiplication satisfies $(ab)c = \sum (x^1\cdot a)[(x^2\cdot b)(x^3\cdot c)]$, for all $a,b,c \in
A$, where $\phi^{-1}=x^1\otimes x^2 \otimes x^3\in H\otimes H \otimes H$ is the inverse of the Drinfeld reassociator of $H$. We say that $H$ acts finitely on $A$ if $\mathrm{Im}(H\rightarrow \End{k}{A})$ is finite dimensional. By the proof of Theorem
\ref{FiniteDimensional}, substituting $A$ by $L(A)$ we get the following

\begin{corollary} Let $A$ be a left $H$-module algebra over a quasi-Hopf algebra $H$ which acts finitely on it. If $A$ has finite left uniform
dimension and is a simple left $A\# H$-module, then it is finite dimensional over $A^H$. \end{corollary}

This applies in particular to finite quasi-Hopf action on non-associative division rings, which are now seen to be finite extensions of their
(associative) subring of invariants. Quasi-Hopf actions on non-associative division rings were considered for example by Albuquerque and Majid in
\cite{AlbuquerqueMajid}.

\section{Critically Compressible Actions}

The compressible condition on $A$ in Theorem \ref{prop5}(2) is stronger than the monoform condition in \ref{Proposition1}(1). In the sequel
of this section, we  will examine how far they actually differ.  This implies that we have to consider some more module theory. The following notions will be need:
a left $R$-module $M$ is called  \textbf{retractable} if $Hom_R (M,N)\neq 0$ for any non-zero submodule $N$ of $M$, while $M$ is called a
\textbf{fully retractable} left $R$-module if for any $g:N\rightarrow M$, $Hom_R (M,N)g\neq 0$ for any non-zero submodule $N$ of $M$.

\begin{proposition}\label{prop1}
Let $M$ be a left $R$-module. 
\begin{enumerate}
	\item[1.] If $M$ is retractable the following statements are equivalent:
	
\begin{enumerate}
	\item[(a)] $M$ is uniform and $\End{R}{M}$ is a domain;
	\item[(b)] $M$ is uniform and any non-zero endomorphism is injective;
	\item[(c)] $\End{R}{M}$ is a left Ore domain. 
\end{enumerate}

\item[2.] Any retractable module $M$ with $\End{R}{M}$ being a domain is compressible.
 
 	\item[3.] The following properties are equivalent:
\begin{enumerate}
	\item[(a)] $M$ is critically compressible;
	\item[(b)] $M$ is compressible and monoform;
	\item[(c)] $M$ is fully retractable and $\End{R}{M}$ is a left Ore domain;
	\item[(d)] $M$ is retractable, $\End{R}{M}$ is a left Ore domain and $\End{R}{\widehat{M}}$ is isomorphic to the division ring of fractions of $\End{R}{M}$.
\end{enumerate}
\end{enumerate}
\end{proposition}

\begin{proof}
(1.) $(a)\Leftrightarrow (b)$ follows from \cite[Theorem 1.4]{RodriguesSant}.
$(a)\Leftrightarrow (c)$ Suppose $0 \neq f,g \in \End{R}{M}$. $M$ is a uniform module and by definition the intersection of any two non-zero
submodules of $M$ have non-zero intersection which in particular results in $N:=Im(f) \cap Im(g) \neq 0$. Since $M$
is retractable there exists an endomorphism $0 \neq h:M \rightarrow N$. Moreover, $\End{R}{M}$ is a domain which implies that any
endomorphism belonging to $\End{R}{M}$ is injective. In particular $h,f$ and $g$ are monomorphisms. Let $\alpha: M\rightarrow (N)f^{-1}$ such that
$(m)\alpha=m'$ if and only if $(m)h=(m')f$. As $f$ is injective, $\alpha$ is well-defined. Moreover, for all $m \in M$ such that $(m)h \in Im(f)$ there
exists an unique $m' \in M$ such that $(m')f=(m)h$. The following composition 
$$M\stackrel{\alpha}{\rightarrow} (N)f^{-1} \stackrel{f}{\rightarrow} N$$
satisfies $(m)\alpha f =(m')f=(m)h$, $\forall m \in M$.
In an analogous way there exists a map $\beta: M \rightarrow (N)g^{-1}$ such that $(m)\beta=m'$ if and only if $(m)h=(m')g$. The composition
$$M\stackrel{\beta}{\rightarrow} (N)g^{-1} \stackrel{g}{\rightarrow} N$$ satisfy $(m)\beta g =(m')g=(m)h$, $\forall m \in M$.
For non-zero maps $f,g \in \End{R}{M}$ we proved the existence of non-zero endomorphisms $\alpha, \beta \in \End{R}{M}$ such that the diagram
$$\xymatrix{{ M}\ar@{-->}[d]_{\beta}\ar@{-->}[r]^{\alpha} \ar[dr]|{h} &
{(N)f^{-1}} \ar^{f}[d] \\
{(N)g^{-1}}\ar_-{g}[r] & N}$$
commutes, i.e. $\alpha f  = \beta g$ and hence $\End{R}{M}$ is a left Ore domain. Conversely assume that $\End{R}{M}$ is a left Ore domain and let $N,L$ be a non-zero submodules of $M$ such
that $N
\cap L=0$. Then $Hom(M,N) \cap Hom(M,L)=Hom(M, N \cap L)=0$. Since $Hom(M,N)$ and $Hom(M,L)$ are left ideals of $\End{R}{M}$ and since $\End{R}{M}$ is
uniform as left $\End{R}{M}$-module we have $Hom(M,N)=0$ or $Hom(M,L)$=0. But $M$ is retractable which implies that $N=0$ or $L=0$. Hence $M$ is
uniform.

(2.) Suppose that $M$ is retractable and $\End{R}{M}$ is a domain. Then for any non-zero submodule $N\subseteq M$, there exists 
$0\neq f:M\rightarrow N$. Since $\Hom{R}{M}{\Ker{f}}f=0$ and $\End{R}{M}$ a domain, $\Hom{R}{M}{\Ker{f}}=0$. Hence $\Ker{f}=0$ as $M$ is retractable.

(3.) $(a)\Leftrightarrow (b)$ follows from \cite[Proposition 1.1]{Zelmanowitz_WeaklyPrimitive}, while $(a)\Leftrightarrow (c)$ follows from
\cite[Proposition 2.4]{RodriguesSant} and (1).
$(a)\Rightarrow (d)$ follows from \cite[11.5(2)]{Zelmanowitz_WeaklyPrimitive}, while $(d)\Rightarrow (a)$ holds since if $f:N\rightarrow M$ is any non-zero endomorphism, then $f$ can be extended to an endomorphism $g$ of $\widehat{M}$ which is injective. Hence also $f$ is injective. \end{proof}

A ring $R$ is a compressible left $R$-module if and only if $R$ is a domain, while $R$ is a critically compressible left $R$-module if and only if it
is a left Ore domain. Hence the free algebra $K \left\langle x,y \right\rangle$ in two indeterminantes over a field $K$ is an example of a
compressible module which is not critically compressible.

Zelmanowitz claimed in \cite{Zelmanowitz_GeneralDensity} that any retractable uniform module whose non-zero endomorphisms are injective would be
critically compressible, but confirmed in \cite{Zelmanowitz_WeaklyPrimitive} that he had no proof for this claim. 
As far as the authors know this claim has not been confirmed up to today. Some partial results on this question have been obtained in \cite{Jeong} and
\cite{RodriguesSant}. However we see from Proposition \ref{prop1} using Rodrigues and Sant'Ana's result \cite{RodriguesSant} that a retractable module
$M$ is uniform with all non-zero endomorphisms injective if and only if $\End{R}{M}$ is a left Ore domain, while $M$ being critically compressible is
equivalent to $M$ being fully retractable and $\End{R}{M}$ being a left Ore domain. {\it Hence Zelmanowitz' claim is equivalent to the question whether a
retractable module $M$ with $\End{R}{M}$ being a left Ore domain is actually fully retractable.} 

The following Lemma will be needed:
\begin{lemma}[Garcia Hernandes - Gomes Pardo] \label{fullyretractable}
Let $M$ be a left $R$-module.
\begin{enumerate}
\item[1.] $M$ is fully retractable if and only if $Hom(N/Tr(M,N),M)=0$, for all $N\subseteq M$;
\item[2.] A retractable self-projective module is fully retractable. 
\end{enumerate}
\end{lemma}
\begin{proof}
{\cite[Proposition 1.1 and 1.2]{GarciaPardo}}
\end{proof}

\begin{lemma} \label{lemself-projective}
A retractable module $M$ with $\End{R}{M}$ being a left Ore domain is critically compressible if $M$ is self-projective or if $M$ is a self-generator.
\end{lemma}
\begin{proof}
Let $M$ be a self-projective module. From Lemma \ref{fullyretractable} it follows that a self-projective retractable module is fully retractable. Since $\End{R}{M}$ is a left Ore domain, we might use Proposition \ref{prop1} $3.(c)\Rightarrow (a)$ to conclude that $M$ is critically compressible.

On the other hand, if $M$ is a self-generator module then $N=Tr(M,N)$, for all $N\subseteq M$. Thus $Hom(N/Tr(M,N),M)=0$ and by Lemma
\ref{fullyretractable}  it follows that $M$ is fully retractable. Hence, from Proposition \ref{prop1} $3.(c)\Rightarrow (a)$ it follows that 
$M$ is critically compressible.
\end{proof}


\subsection{Compressible operator actions}
Let us suppose again that we have an intermediate algebra $L(A) \subseteq B \subseteq End(A)$, where $A$ is a $k$-algebra. Using the identification
$\Hom{B}{A}{-} \cong (-)^B$ we see that $A$ is a retractable left $B$-module if and only if $A^B\cap I \neq 0$ for any non-zero $B$-stable left ideal
$I$ of $A$. In the literature this property is also referred to as $A^B$ being large in $A$. Thus we reencounter a property that is for example known
for semiprime PI-algebras $A$ which have a large center, i.e. $A$ is a retractable $A^e$-module. Another instance is Bergman and Isaacs theorem on finite group action which says that if a finite group $G$ act on a semiprime ring $A$ which is $|G|$-torsionfree, then $A^G$ is large in $A$, i.e. $A$ is a retractable left $A*G$-module. 

Recall that a multiplicatively closed subset $S$ of right non-zero divisors of a ring $R$ is called {\it left permutable} if $Sa \cap Rs$ is not empty, for every $a \in R$ and $s \in S$.

Additional to the module theoretic conditions that imply $_BA$ to be critically compressible we have the following:

\begin{lemma}\label{Lemma_OreCondition}
Let $L(A) \subseteq B \subseteq End (A)$. Then $A$ is a critically compressible left $B$-module if $A^B$ is a domain and large in $A$ such that 
\begin{enumerate}
\item $udim(A)=1$ or 
\item $A^B\setminus\{0\}$ is a left  permutable set in $A$ and consists of left non-zero divisors in $A$.
\end{enumerate} 
\end{lemma}

\begin{proof}
(1) In order to prove that $A$ is a critically compressible $B$-module it is enough to show that any partial endomorphism of $A$ is injective. Let
$f:I\rightarrow A$ be a non-zero $B$-linear map with $0\neq {_BI} \subseteq {_BA}$ and denote $K:=ker(f)$ and $J:=(I)f$. 
 Now suppose $K \neq 0$.  Since $_BA$ is retractable there exists an element $0 \neq x \in A^B \cap K$ and an element $0 \neq y \in A^B \cap J$ which implies that there exists $i \in I$ such that $y=(i)f$. As $A$ is uniform, $Ax \cap Ai \neq 0$ and hence there exists $z,t \in A$ such that $zx=ti \neq 0$.
But $ty=t(i)f=(ti)f=(zx)f=z(x)f=0$.
Since $0=ty=(t)R_y$ and $R_y$ is a $B$-linear endomorphism of $A$, we have by Proposition \ref{prop1} that $R_y$ is injective and hence $t=0$. Thus
$zx=ti=0$ which is a contradiction. Thus $K=0$ and we conclude that $f$ is injective.

(2) Let $f:I\rightarrow A$ be a $B$-linear map with $0\neq {_BI} \subseteq {_BA}$ and denote $K:=ker(f)$. In order to prove that $A$ is a critically
compressible $B$ -module it is enough to show that $f$ is a monomorphism. Since $_BA$ is retractable then for all $0 \neq L \subseteq A$, $A^B \cap L
\neq 0$. In particular if $K$ is a non-zero submodule of $A$, there exists an element $0 \neq x \in A^B \cap K$. Since $A^B\setminus\{0\}$ is left permutable set in $A$, we have that for all $a \in I$ there exists $y \in A^B \setminus\{0\}$ and $b \in A$ such
that $bx=ya$. But $0=b(x)f=(bx)f=(ya)f=y(a)f$. By the hypthesis $A^B\setminus\{0\}$ consists of left non-zero divisor set in $A$. Thus $(a)f=0$ or
$K=0$. Hence we conclude that $f$ is injective.
\end{proof}

As an application we get  a characterisation of the critically compressibility of $A$ for central invariants $A^B\subseteq Z(A)$.

\begin{corollary} \label{prop3}
Let $L(A) \subseteq B \subseteq End (A)$ such that $A^B\subseteq Z(A)$. Then  $A$ is a critically compressible $B$-module if and only if $A^B$ is an
integral domain and large in $A$.
\end{corollary}

\begin{proof} By Proposition \ref{prop1} the elements of $A^B$ are non-zero divisors on the right if $_BA$ is retractable and $A^B\simeq \End{B}{A}$
is a domain. Since $A^B\subseteq Z(A)$, the non-zero elements of $A^B$ are left and right non-zero divisors and form a left permutable set in $A$. By
Lemma \ref{Lemma_OreCondition}, $A$ is critically compressible.
The converse is clear by \ref{prop1}.
\end{proof}

\subsection{Group Actions and Group gradings}

In this section we will apply our earlier results to group actions and group gradings.

\begin{proposition}
Let $G$ be a finite group acting as automorphism on $A$ such that $A$ is $|G|$-torsionfree. Then $A$ is a critically compressible left $A\ast G$-module if and only if $A^G$ is a left Ore domain and large in $A$. In this case $udim(A)\leq |G|$ and $dim({\hat{A}}_{\Delta}) \leq |G|^2$  where $\widehat{A}$ is the self-injective hull of $A$ and $\Delta$ is isomorphic to the division ring of fractions of $A^G$. 
\end{proposition}

\begin{proof}
We may assume that $|G|$ is invertible in $A$, otherwise we localize $A$ by the powers of $|G|$, obtaining an algebra in which $|G|$ is invertible.
If $|G|^{-1}\in A$ then $A$ a self-projective $A\ast G$-module (see \cite{Wisbauer}). 
Since $A$ is a retractable left $A\ast G$-module and $A^G \cong End_{A \ast G}(A)$ is a left Ore domain, 
it follows from Lemma \ref{lemself-projective} that $A$ is a critically
compressible left $A \ast G$-module. Moreover, from Proposition \ref{prop1} $1.c) \Rightarrow b)$ we know that $A$ is an uniform $A \ast G$-module which
implies that $udim(_{A \ast G}A)=1$.  By {\cite[Lemma 1.3]{Passman}} we have $udim(A)\leq |G| \cdot udim({_{A \ast G}A})$ and hence $udim(A)\leq |G|$. Now using this fact, the relation $dim({\hat{A}}_{\Delta}) \leq |G|^2$ follows directly from Theorem \ref{prop5}.

\end{proof}

Let $G$ be a group. A $G$-graded $k$-algebra $A$ is a $k$-algebra with a decomposition 
$A=\bigoplus_{g\in G} A_g$ into subspaces $A_g$, $g\in G$, such that $A_gA_h\subseteq A_{gh}$, for all $g,h \in G$. 
For each $g\in G$ let $\pi_g:A\rightarrow A_g$ be the projection onto the $g$-component. 
Define $B=\langle L(A)\cup \{\pi_g \mid g\in G\}\rangle \subseteq \End{k}{A}$. 
The left $B$-submodules are precisely the $G$-graded left ideals of $A$ and if $e$ denotes the neutral element of $G$, then $A_e = A^B$ is a subring
of $A$ isomorphic to $\End{B}{A}$. 

\begin{proposition}
Let $G$ be a group and $A$ be a $G$-graded $k$-algebra. Then $A$ is a critically compressible left $B$-module if and only if $A_e$ is a left Ore
domain and large in $A$.
\end{proposition}

\begin{proof} 
We start by showing that a $G$-graded $k$-algebra $A$ is a self-projective $B$-module. In order to prove that we will show that for any $G$-graded left ideal $I$ of $A$ we have $Hom_B(A,A/I)=End_B(A) p_I$, where $p_I:A \rightarrow A/I $ is the projection map. Since  
$$(A_e)p_I = \left( A_e +I \right) / I \cong A_e / \left( I\cap A_e\right) = A_e/I_e \cong \left( A/I \right)_e = (A/I)^B$$
and $End_B(A)\cong A_e$, we have $$End_B(A) p_I \cong (A_e)p_I = (A/I)^B \cong Hom_B(A,A/I).$$
Since $A$ is a self-projective $B$-module it follows from Lemma \ref{lemself-projective} that $A$ is a critically compressible module.
The converse follows from \ref{prop1}(3) and $\End{B}{A}\simeq A_e$.
\end{proof}

\subsection{Lie actions}

A derivation $\delta$ of $A$ is called \textbf{locally nilpotent} if for any $a \in A$ there exists $n>0$ such that $\delta^n(a)=0$. Recall $A^\delta=\Ker{\delta}$. For each non-zero element $a\in A$ one defines then its degree (with
respect to the derivation $\delta$) as $deg(a)=\mathrm{min}\{n\in \NN_0: \delta^{n+1}(a)=0\}$.

\begin{lemma}\label{prop_retract_derivations}
Let $A$ be a $k$-algebra  and let $\delta$ be a locally nilpotent derivation of $A$. Then $A$ is a retractable $A[z, \delta]$-module.
\end{lemma}

\begin{proof} By the correspondence $\Hom{A[z:\delta]}{A}{ - } \simeq ( - )^\delta$ we have to show that  every  non-zero $\delta$-stable left ideal
$I$ of $A$ contains a non-zero constant element, i.e. $I\cap A^\delta \neq 0$. Since $\delta$ is locally nilpotent, for any $0\neq a \in I$ with
$n=deg(a)$, we have $0\neq \delta^n(a)\in A^\delta\cap I$.  
\end{proof}

We will show that if $A$ is an algebra over a field $k$ with a locally nilpotent derivation $\delta$ having an element $x\in A$ such that $\delta(x)=1$,
then $A$ is a self-projective $A[z;\delta]$-module if $char(k)=0$ or $\delta^p=0$ in case $char(k)=p>0$. In order to do so we first show that any
element of $A$ can be written as polynomial in $x$ with coefficients in $A^B$ (see also \cite{Bavula}).

\begin{lemma} \label{cor2}
Let $A$ be an algebra over a field $k$ and $\delta$ a locally nilpotent derivation of $A$ or $\delta^p=0$ in case of $char(k)=p$ such that
$\delta(x)=1$, for some $x$ in $A$. Then any element $a \in A$ of degree $n$ can be written as $a=\sum^{n}_{i=0} c_i x^i$, for some $c_i \in A^\delta$.
\end{lemma}

\begin{proof}
Let $0\neq a\in A$ and let $n=deg(a)$. Note that if $char(k)=p$ and $\delta^p=0$ then $n<p$. 
We will use induction on $n$: If $n=0$ then $a \in A^\delta$ and there is nothing to prove.
Suppose that $n>0$. We claim that $deg(a-\frac{1}{n!} \delta^n(a)x^n)<n$. To see this note that $\delta^n\left( a-\frac{1}{n!} \delta^n(a)x^n
\right) = \delta^n(a)-\frac{1}{n!} \delta^n(\delta^n(a)x^n)$.
Since $deg(a)=n$,  $\delta(\delta^n(a))=\delta^{n+1}(a)=0$ which implies that $\delta^n(a) \in A^{\delta}$. Moreover, for any $y \in A$ and $c \in
A^{\delta}$ we have $\delta(cy)=c \delta(y)$. In particular,
$$\delta^n\left( a-\frac{1}{n!} \delta^n(a)x^n \right) =  \delta^n(a)-\frac{1}{n!} \delta^n(a)\delta^n(x^n) = \delta^n(a)-\frac{1}{n!}\delta^n(a)
n! = 0. $$
Hence $deg(a-\frac{1}{n!} \delta^n(a)x^n)<n$. Now suppose that for all $a' \in A$ such that $m=deg(a')<n$, there exist $c_0, c_1, \ldots, c_m \in
A^\delta$ such that $a'=\sum^{m}_{i=0} c_i x^i$. As $deg(a-\frac{1}{n!} \delta^n(a)x^n)<n$,  by the induction hypothesis there are elements
$c_0,\ldots, c_{n-1} \in A^\delta$ with  $a-\frac{1}{n!} \delta^n(a)x^n =\sum^{m}_{i=0} c_i x^i$, equivalently $ a =\frac{1}{n!}
\delta^n(a)x^n + \sum^{m}_{i=0} c_i x^i$, where $\frac{1}{n!} \delta^n(a)\in A^\delta$.
\end{proof}

\begin{proposition}\label{prop6}
Let $A$ be an algebra over a field $k$ and $\delta$ a locally nilpotent derivation of $A$ or $\delta^p=0$ in case of $char(k)=p$ such that
$\delta(x)=1$, for some $x$ in $A$. Then $A$ is a self-projective $A[z,\delta]$-module.
\end{proposition}

\begin{proof}
We have to prove that $A$ is a self projective $A[z,\delta]$-module which is equivalent to prove that $(A/I)^\delta=\frac{A^\delta+ I}{I}$, for all
$\delta$-stable left ideals $I$ of $A$. Thus we need to show that 
$$\forall a \in A: a+I \in (A/I)^\delta \Rightarrow \exists y \in A^\delta : a+I=y+I.$$
Since ${(A/I)}^\delta=\left\{a+I:\delta(a)+I=0+I \right\}=\left\{a+I: \delta(a) \in I \right\}$ the above condition can be rewritten as 
\begin{equation}\label{self-proj-lie}\forall a \in A: \delta(a) \in I \Rightarrow \exists y \in A^\delta : a-y \in I.\end{equation}
We wil prove this by induction on the degree of $a$. Let $a \in A$ be such that $\delta(a)\in I$. If $deg(a)=0$ then $a \in A^\delta$ and thus $y=a$ is
a solution for (\ref{self-proj-lie}). Suppose that $deg(a)=n>0$ and that for all $b\in A$ such that $deg(b)<n$ and $\delta(b)\in I$, there exists a solution $y\in A^\delta$
such that $b-y \in I$. By Lemma \ref{cor2}, $a=\sum^{n}_{i=0} c_i x^i$ with $c_i \in A^\delta$ and $n=deg(a)$. The derivation of $a$ is given
by $\delta(a)=\sum^{n}_{i=0} c_i \delta(x^i)=\sum^{n}_{i=0} c_i ix^{i-1}$. In particular, if we consider $b=a-\frac{1}{n}x\delta(a)$ it follows that
\begin{eqnarray*}
\delta(b)&=&\delta(a)-\frac{1}{n}\delta(x\delta(a))\\
&=&\delta(a)-\frac{1}{n}(\delta(x) \delta(a)+x \delta^2(a))\\
&=&\delta(a)-\frac{1}{n} \delta(a)-\frac{1}{n}x \delta^2(a)\\
&=& (1-\frac{1}{n})\delta(a)-\frac{1}{n}x \delta^2(a)
\end{eqnarray*}
 which belongs to $I$ since $\delta(a)\in I$ and $I$ is $\delta$-stable. In order  to show that $deg(b)<n$ we will write $b$ using
the above equalities for $a$ and $\delta(a)$. Denote the commutator of two elements $r,s \in A$ by $[r,s]=rs-sr$. Thus 
\begin{eqnarray*}
b &=&\sum^{n}_{i=0} c_i x^i-\frac{1}{n}x\left(\sum^{n}_{i=1} c_i ix^{i-1}\right) \\
& =& \sum^{n-1}_{i=0} c_i x^i+c_n x^n- \sum^{n-1}_{i=1} \left(xc_i \frac{i}{n}x^{i-1}\right)- x
c_n x^{n-1}\\
& =& \sum^{n-1}_{i=0} c_i x^i+c_n x^n - \sum^{n-1}_{i=1} \left(c_i x + [x,c_i]\right)\frac{i}{n}x^{i-1} - c_n x x^{n-1}- [x,c_n]x^{n-1}\\ 
& =& \sum^{n-1}_{i=0} c_i x^i -\sum^{n-1}_{i=1} c_i \frac{i}{n} x^i - \sum^{n-1}_{i=1} [x,c_i]\frac{i}{n}x^{i-1} - [x,c_n]x^{n-1}\\
& =& \sum^{n-1}_{i=0} c_i \left(1-\frac{i}{n}\right) x^i - \sum^{n}_{i=1} \frac{i}{n} [x,c_i]x^{i-1}.
\end{eqnarray*}
Hence $deg(b)<n$ since $[x,c_i]\in A^\delta$ for all $c_i\in A^\delta$. By induction, there exist $y \in A^\delta$ such that $b-y \in I \Leftrightarrow a- \frac{1}{n} x \delta(a) -y \in I$. Since $\frac{1}{n} x \delta(a) \in I$, $a-y \in I$, as required. 
\end{proof}

\begin{corollary}\label{Corollary_Lie}
Let $A$ be an algebra over a field $k$ and $\delta$ a locally nilpotent derivation of $A$ or $\delta^p=0$ in case of $char(k)=p$ and suppose that there exists $x\in A$ such that $\delta(x)=1$. Then $A^\delta$ is a left Ore domain if and only if $A$ is a critically compressible $A[z,\delta]$-module.
\end{corollary}

\begin{proof}
Suppose that $A^\delta$ is a left Ore domain. From Proposition \ref{prop6}, $A$ is a self-projective $A[z,\delta]$-module. Since $\delta$ is
locally nilpotent, $A$ is a rectratable $A[z,\delta]$-module. From Proposition \ref{fullyretractable} $2.$ it follows that
$A$ is fully retractable as $A[z,\delta]$-module. Hence, by Proposition \ref{prop1} $3.$, $A$ is a critically compressible $A[z,\delta]$-module if
and only if $A^\delta$ is a left Ore domain.
\end{proof}

Actually $A$ is a domain in case $char(k)=0$ and $_{A[z,\delta]}A$ is critically compressible. Hence in this case condition (2) of Theorem \ref{Proposition1} implies condition (1).

\begin{proposition}
Let $A$ be an algebra over a field $k$ of characteristic zero and $\delta$ be a locally nilpotent derivation of $A$ such that $\delta(x)=1$, for some
$x$ in $A$. Then the following statements are equivalent:
\begin{enumerate}
\item[(a)] $A^\delta$ is a left Ore domain;
\item[(b)] $A$ is a left Ore domain;
\item[(c)] $A$ is a critically compressible $A[z,\delta]$-module.
\end{enumerate}
\end{proposition}

\begin{proof}
$(a)\Leftrightarrow (c)$ follows from Corollary \ref{Corollary_Lie}.

$(a)\Rightarrow (b)$ 
Let $A^\delta$ be a left Ore domain. By Lemma \ref{lemma1}, $A^\delta$ is a monoform left $A^\delta$-module. Now from {\cite[Lemma 2.1.]{Bavula}}
it follows that $A=A^\delta\left[x,d\right]$ is the Ore extension with coefficients from the algebra $A^\delta$ and where the derivation $d$ of the
algebra $A^\delta$ is the restriction of the inner derivation of the algebra $A$ to its subalgebra $A^\delta$ defined by $d(a):=\left[x,a\right]$, for all $a\in A$. Thus by {\cite[Proposition 1.3.]{BellGoodearl}} we have that $A=A^\delta[x,d]$ is a monoform left $A$-module. Finally, by Lemma \ref{lemma1} $A$ is a left Ore domain.

$(b) \Rightarrow (c)$ Let $A$ be a left Ore domain, then by Lemma \ref{lemma1}, $A$ is a monoform $A[z,\delta]$-module. In particular 
$A$ is a uniform $A[z,\delta]$-module. By \ref{prop_retract_derivations},  $A$ is a retractable $A[z,\delta]$-module and hence by Proposition \ref{prop1} $\End{A[z,\delta]}{A}\simeq A^\delta$ is a left Ore domain.
\end{proof}

However in positiv characteristic $A$ might not be even reduced.
\begin{example}
Let $A=k\left[x \right]$ be the commutative polynomial ring in one variable, where k is a field of characteristic $p>0$. Consider the
locally nilpotent derivation $\delta$ of $A$ defined by $\delta=\frac{\partial}{\partial x}$. Now let $I=\left\langle x^p \right\rangle$. Since
$char(k)=p$ we have $\delta(x^p)=px^{p-1}=0$. Moreover, for all $f(x) \in A$, $\delta (f(x)\cdot x^p)=\delta(f(x)) \cdot x^p \in I$. Hence $I$ is $\delta$-stable and we can define a derivation $\bar{\delta}$ in $\bar{A}=A/I$ such that
$\bar{\delta}(\bar{f(x)})=\delta(f(x))+I$. In particular, for $\bar{x}$ we have $\bar{\delta}(\bar{x})=\delta(x)+I=\bar{1}$. On the other hand 
${\bar{\delta}}^{p}(f)=0$, for all $\bar{f(x)} \in \bar{A}$. Since ${\bar{A}}^{\bar{\delta}}=k+I\cong k$ is a field, $\bar{A}$ is a simple left $\bar{A}[z;\delta]$-module by Lemma \ref{prop_retract_derivations}, but $\bar{A}$ is not a domain (not even reduced).
\end{example}

\end{document}